\documentclass[12pt]{amsart} 
\usepackage{verbatim}
\usepackage{epsfig}

\scrollmode

\def\AAP{Ann.\ Appl.\ Prob.\ }
\def\FOCM{Found.\ Comput.\ Math.\ }
\def\JAMSA{J.\ Appl.\ Math.\ Stochastic Anal.\ } 
\def\JC{J.\ Complexity\ }
\def\JCAM{J.\ Comput.\ Appl.\ Math.\ }
\def\MC{Math.\ Comp.\ }
\def\PA{Potential Analysis\ }%
\def\SJNA{SIAM J.\ Numer.\ Anal.\ }
\def\SPA{Stochastic Processes Appl.\ }
\def\SSR{Stochastics Stochastics Rep.\ }
\newcommand{\R}{{\mathbb R}}  
\newcommand{\N}{{\mathbb N}} 
\newcommand{\spann}{\operatorname{span}}
\newcommand{\scp}[2]{\langle #1 , #2 \rangle}
\newcommand{\lscp}[2]{\bigl\langle #1 , #2 \bigr\rangle}
\newcommand{\LL}{{\mathcal L}} 
\newcommand{\J}{{\mathcal J}}
\newcommand{\I}{{\mathcal I}}
\newcommand{\K}{{\mathcal K}}
\newcommand{\bj}{{\boldsymbol{j}}}
\newcommand{\bi}{{\boldsymbol{i}}}
\newcommand{\Xh}{\widehat{X}}
\newcommand{\Yh}{\widehat{Y}}
\newcommand{\Xb}{\overline{X}}
\newcommand{\Yb}{\overline{Y}}
\newcommand{\Xt}{\widetilde{X}}
\newcommand{\Yt}{\widetilde{Y}}
\newcommand{\XX}{\mathfrak{X}}
\newcommand{\at}{a}
\newcommand{\ah}{a}
\newcommand{\Bb}{\overline{B}}
\newcommand{\pb}{\overline{\psi}}
\newcommand{\ph}{\widehat{\psi}}
\renewcommand{\epsilon}{\varepsilon}

\theoremstyle{plain}
\newtheorem{theorem}{Theorem}
\newtheorem{prop}{Proposition}
\newtheorem{lemma}{Lemma}

\theoremstyle{definition}
\newtheorem{rem}{Remark}

\begin{document}

\title[An Implicit Euler Scheme with Non-uniform Time Discretization]%
{An Implicit Euler Scheme with Non-uniform Time Discretization for 
Heat Equations with \\ Multiplicative Noise}

\author[]
{Thomas M\"uller-Gronbach}
\address{Institut f\"ur Mathematische Stochastik\\
Fakult\"at f\"ur Mathematik\\
Universit\"at Magdeburg\\
Postfach 4120\\
39016 Magdeburg\\
Germany}
\email{gronbach@mail.math.uni-magdeburg.de}

\author[]
{Klaus Ritter}
\address{Fachbereich Mathematik\\
Technische Universit\"at Darmstadt\\
Schlo\ss gartenstra\ss e 7\\
64289 Darmstadt\\
Germany}
\email{ritter@mathematik.tu-darmstadt.de}

\date{April 2006}

\begin{abstract}
We present an algorithm for solving stochastic
heat equations, whose key ingredient is a non-uniform
time discretization of the driving Brownian motion $W$.
For this algorithm we derive an error bound
in terms of its number of evaluations of one-dimensional
components of $W$. The rate of convergence depends on
the spatial dimension of the heat equation and on the decay of the 
eigenfunctions of the covariance of $W$.
According to known lower bounds, our algorithm
is optimal, up to a constant, and this optimality 
cannot be achieved by uniform time discretizations.
\end{abstract}

\maketitle

\section{Introduction}

A common technique for the numerical solution of stochastic evolution 
equations is an It\^{o}-Galerkin approximation,
which turns the corresponding infinite-dimensional system of stochastic
differential equations (SDEs) into a finite-dimensional one. The
latter is then discretized in time and 
approximately solved by, e.g., an Euler scheme.

More generally, every numerical algorithm for an evolution equation
eventually has to discretize the driving stochastic process,
which frequently is assumed to be a Brownian motion on an
infinite-dimensional Hilbert space, 
in space and time. The vast majority of algorithms
for stochastic evolution equations as well as for SDEs apply a uniform 
time discretization. This means that a finite number of
one-dimensional components of the driving process are
evaluated (simulated) at time instances $\ell/n$ with a common
step-size $1/n$.

In this paper we present and analyze a non-uniform time discretization
for a stochastic
heat equation
\begin{equation}\label{g1}
\begin{aligned}
dX(t) &= \Delta X(t) \, dt + B(X(t))\, dW(t),\\
X(0) &= \xi
\end{aligned}
\end{equation}
on the Hilbert space $H = L_2(\left]0,1\right[^d)$.
As a key assumption, the system $(h_\bi)_{\bi \in \N^d}$
of eigenfunctions of the trace class
covariance $Q$ of the Brownian motion $W$ coincides with the
system of eigenfunctions of the Laplace operator $\Delta$.
A finite number of scalar Brownian motions $\scp{W}{h_\bi}$ is selected,
and each of them is evaluated with step-size $1/n_\bi$ depending on
its variance. Based on these data, a properly defined implicit Euler scheme
is employed to compute an approximation $\Xh^*_N$
to a finite number of components
$\scp{X}{h_\bj}$ of the solution $X$. 
Here $N$ denotes the total number of evaluations of scalar Brownian 
motions $\scp{W}{h_\bi}$ used by $\Xh^*_N$, up to a constant.

Our main result is an upper bound for the error
\[
e(\Xh^*_N) = \left(E \left(
 \int_0^1 \| X(t) - \Xh^*_N(t) \|^2 \, dt \right) \right)^{1/2}
\]
of $\Xh^*_N$ in terms of $N$. The rate of convergence
depends on the spatial dimension $d$ and on the decay 
of the eigenvalues of the covariance $Q$.
Assume, for simplicity, that
\[
Q h_\bi = |\bi|_2^{-\gamma} \cdot h_\bi
\]
for some
$\gamma \in \left] d, \infty \right[ \setminus \{ 2d\}$, and put
\[
\alpha^*(\gamma,d) = 
\frac{1}{2}-\frac{(2d-\gamma)_+}{2(d+2)}.
\]
Then
\[
e(\Xh^*_N) \leq c_1 \cdot N^{-\alpha^*(\gamma,d)}
\]
with a constant 
$c_1 > 0$ that only depends on $d$, $\gamma$, $B$, and $\xi$, 
see Theorem \ref{t2}.

Actually, this upper bound is best possible, not only for
the specific algorithm $\Xh_N^*$ but for any algorithm that
uses at most a total of $N$ evaluations of the scalar
Brownian motions $\scp{W}{h_\bi}$: there exists a constant
$c_2 > 0$ that only depends on $d$, $\gamma$, $B$, and $\xi$ such that
\[
e(\Xh_N) \geq c_2 \cdot N^{-\alpha^*(\gamma,d)}
\]
for any such algorithm. In general,  one cannot
achieve the optimal rate $\alpha^*(\gamma,d)$ by 
any sequence of algorithms that use a uniform discretization.
See M\"uller-Gronbach, Ritter (2006).

In the context of stochastic partial differential
equations, implicit (Euler) schemes based on uniform
time discretizations are studied, e.g., by
Gy\"ongy (1999),
Kloeden, Shott (2001),
Hausenblas (2002, 2003),
Millet, Morien (2005),
Walsh (2005),
and
Yan (2005).
Non-uniform time discretizations are studied for the
first time by M\"uller-Gronbach, Ritter (2006). 
In the latter paper, non-uniform time discretizations are used
for the numerical solution of heat equations with additive noise, 
i.e.,
$B$ is a function of the time $t$ but not of the current value $X(t)$
of the evolution. In this case
the solution $X$ is a Gaussian process and conditional
expectations become feasible as a computational tool.
This is no longer true for equations with multiplicative noise,
as studied in the present paper. Instead, the algorithm introduced
in the present paper is a general-purpose algorithm.

Optimality results, as stated here for the algorithm $\Xh_N^*$, 
require lower bounds that are valid for all (or at least
a broad class) of algorithms. For stochastic evolution
equations the first such lower bound is due to 
Davie, Gaines (2001), who consider a particular case
of \eqref{g1} in spatial dimension $d=1$
with a space-time white noise. See
M\"uller-Gronbach, Ritter (2006) for lower bounds for
equations \eqref{g1} in general, with space-time
white noise for $d = 1$ and trace class noise for $d \geq 1$.

Our results show the principal significance of non-uniform time
discretizations for the numerical solution of stochastic 
evolution equations. 
Non-uniform and even sequentially computed time-discretizations
are studied, too, for finite-dimensional systems of
SDEs. Here, as a rule, those time-discretizations
do not improve the order of convergence, but only
the asymptotic constants. However, improvements
may be substantial on the level of constants, see Cambanis, Hu (1996), 
Hofmann, M\"uller-Gronbach, Ritter (2001), and 
M\"uller-Gronbach (2002, 2004).

We outline the content of the paper.
In Section \ref{s2} we formulate the assumptions on
the heat equation \eqref{g1} and briefly discuss
existence and uniqueness of a mild solution.
Our algorithm is defined in Section \ref{s3}. Error bounds
and optimality properties are stated in Sections \ref{s4} and
\ref{s5}, resp., and proofs are given in Section \ref{sp}.

\section{Assumptions}\label{s2}

We study stochastic heat equations \eqref{g1}
on the Hilbert space $H = L_2(\left]0,1\right[^d)$.
Here $\xi \in H$ for the initial value,
and $\Delta$ denotes the Laplace 
operator with Dirichlet boundary
conditions on $H$.
Hence $\Delta h_\bi = - \mu_\bi \cdot h_\bi$
with 
\[
h_\bi (u) = 
2^{d/2} \cdot \prod_{\ell=1}^d \sin( i_\ell \pi u_\ell)
\]
and 
\[
\mu_\bi = \pi^2 \cdot |\bi|_2^2,
\]
where $|\bi|_2$ is the Euclidean norm of $\bi \in \N^d$.

Moreover, $W=(W(t))_{t\in [0,1]}$ denotes a Brownian motion on $H$,
whose covariance $Q : H \to H$ is a trace class operator.
Specifically, we assume that $Q h_\bi = \lambda_\bi \cdot h_\bi$ with
\[
\lambda_\bi = \lambda(|\bi|_2)
\]
for some non-increasing and regularly varying function 
\[
\lambda : \left[1,\infty\right[ \to \left]0,\infty\right[
\]
of index $-\gamma$, where 
\[
\gamma \in \left[ d, \infty\right[ \setminus \{ 2d\}
\]
and
\[
\int_1^\infty \lambda(r) \cdot r^{d-1} \, dr < \infty.
\]
Note that the latter property always holds if $\gamma > d$.

Let $\scp{\cdot}{\cdot}$ and $\|\cdot\|$ denote the inner product
and the norm, respectively, in $H$, 
and consider the Hilbert space $H_0=Q^{1/2} (H)$,
equipped with the inner product
$(Q^{1/2} h_1,Q^{1/2} h_2) = \scp{h_1}{h_2}$
for $h_1, h_2 \in H$. 
Furthermore, let $\LL = \LL_2(H_0,H)$
denote the class of Hilbert-Schmidt operators from $H_0$ into $H$, 
equipped with the Hilbert-Schmidt norm $\|\cdot \|_\LL$. 
We assume that the mapping $B : H \to \LL$ is given by pointwise 
multiplication and a Nemytskij operator, i.e.,
\[
B(x) h = T_g (x) \cdot h
\]
for $x \in H$ and $h \in H_0$, where
\[
T_g (x) = g \circ x
\]
with $g \in C^1(\R)$ such that 
\[
\|g^\prime\|_\infty < \infty.
\]

\begin{rem}\label{rem1}
We briefly comment on the existence of a mild solution of equation 
\eqref{g1} under the above conditions on $B$.

Note that
$T_g: H\to H$, see Appell, Zabrejko (1990, Thm.\ 3.1). Furthermore, 
$H_0\subset L_\infty(\left]0,1\right[^d)$, since 
$\sup_{\bi\in\N^d}\|h_\bi\|_\infty<\infty$, see Manthey, Zausinger 
(1999, Lemma 2.2). Hence $B(x) h\in H$ for every $h\in H_0$. Moreover,
$(\lambda_\bi^{1/2}\cdot h_\bi)_{\bi\in\N^d}$ is a complete 
orthonormal system in $H_0$ and
\[
\sum_{\bi\in\N^d} \|B(x)\lambda_\bi^{1/2}h_\bi\|^2 \le 
\sup_{\bi\in\N^d}\|h_\bi\|^2_\infty\cdot \sum_{\bi\in\N^d}\lambda_\bi
\cdot \|T_g(x)\|^2.
\]
Consequently, $B(x)\in\LL$ for every $x\in H$. Clearly,
\begin{equation}\label{g13}
\|T_g (x) - T_g (y)\| \leq \|g^\prime\|_\infty \cdot \|x-y\|
\end{equation}
for $x,y \in H$, which yields 
\begin{equation}\label{e1}
\|B(x)-B(y)\|_\LL \le K\cdot \|x-y\|
\end{equation}
with $K = \sup_{\bi\in\N^d}\|h_\bi\|_\infty\cdot
\sum_{\bi\in\N^d}\lambda_\bi \cdot \|g^\prime\|_\infty$. 
Thus $B: H \to \LL$ is Lipschitz continuous.

Consider the semigroup $(S(t))_{t \geq 0}$ on $H$ that is
generated by $\Delta$, i.e.,
\[
S(t) h_\bi = \exp(- \mu_\bi t) \cdot h_\bi.
\]
{}From \eqref{e1} it follows that there exists a continuous process 
$(X(t))_{t \in [0,1]}$ with values 
in $H$, which is 
adapted to the underlying filtration,
such that, for every $t \in [0,1]$,
\[
X(t) = S(t) \xi 
+ \int_0^t S(t-s) B(X(s)) \, dW(s) 
\]
holds a.s. 
This process is uniquely determined a.s., and it
is called the mild solution of equation \eqref{g1}.
Furthermore,
\begin{equation}\label{g3}
\sup_{t \in [0,1]} E\| X(t)\|^2 \le c_1,
\end{equation}
where the constant $c_1>0$ only depends on $d, \xi, \lambda$ and $g$.
See Da Prato, Zabczyk (1992, Sec.\ 7.1).
\end{rem}

\section{The Algorithm}\label{s3}

We construct and analyze an algorithm 
that is built from the following ingredients:
\begin{itemize}
\item[(i)] an It\^{o}-Galerkin approximation of the stochastic heat
equation,
\item[(ii)] a non-uniform time discretization of the corresponding
      finite-dimensional Brownian motion,
\item[(iii)] a drift-implicit Euler scheme.
\end{itemize}
Put
\[
\beta_\bi(t) = \lambda_\bi^{1/2}\cdot \scp{W(t)}{h_\bi}
\]
for $\bi\in\N^d$ and $t\in [0,1]$. Then $(\beta_\bi)_{\bi\in\N^d}$ is
an independent family of standard one-dimensional Brownian motions.
Let
\[
Y_\bj(t) = \scp{X(t)}{h_\bj}
\]
for $t \in [0,1]$ and $\bj \in \N^d$.
The real-valued processes $Y_\bj = (Y_\bj(t))_{t \in [0,1]}$ 
satisfy the bi-infinite system
\begin{equation*}
\begin{aligned}
d Y_\bj (t) &= - \mu_\bj \, Y_\bj (t) \, dt + 
\sum_{\bi \in \N^d} \lambda_\bi^{1/2} \cdot
\lscp{B(X(t)) h_\bi}{h_\bj} \, d \beta_\bi(t)
\\
Y_\bj (0) &= \scp{\xi}{h_\bj}
\end{aligned}
\end{equation*}
of stochastic differential equations.
For any choice of finite sets $\I, \J \subseteq \N^d$
an It\^o-Galerkin approximation $\Xb = (\Xb(t))_{t \in [0,1]}$ to 
$X$ is given by
\[
\Xb (t) = \sum_{\bj \in \J} \Yb_\bj (t) \cdot h_\bj
\]
with real-valued processes $\Yb_\bj = (\Yb_\bj(t))_{t \in [0,1]}$ 
that solve the finite-dimensional system 
\begin{equation}\label{g16}
\begin{aligned}
d \Yb_\bj (t) &= - \mu_\bj \, \Yb_\bj (t) \, dt + 
\sum_{\bi \in \I} \lambda_\bi^{1/2} \cdot
\lscp{B(\Xb(t)) h_\bi}{h_\bj} \, d \beta_\bi(t) \\
\Yb_\bj (0) &= \scp{\xi}{h_\bj}.
\end{aligned}
\end{equation}

We apply a drift-implicit Euler scheme to the finite-dimensional
system \eqref{g16}. This scheme is based on a non-uniform
discretization of the corresponding finite-dimensional Brownian
motion, since $\beta_\bi$ will be evaluated with step-size
$1/n_\bi$ depending on $\bi \in \I$. Accordingly, put
\begin{equation}\label{disc}
\phantom{\qquad \quad \ell = 0,\dots,n_\bi.}
t_{\ell,\bi} = \ell/n_\bi, 
\qquad \quad \ell = 0,\dots,n_\bi.
\end{equation}
A good choice of the integers $n_\bi \in \N$, together with sets 
$\I$ and $\J$, will be presented in Section \ref{s4}. 

In order to understand the construction of this scheme better we
first consider
a uniform discretization, i.e., $t_\ell = t_{\ell,\bi} = \ell / n$
with a common step-size $1/n$ for all $\bi \in \I$. In this case
the drift-implicit Euler scheme is given by
\[
\Yh_\bj (t_\ell) = \Yh_\bj (t_{\ell-1}) - \mu_\bj\,
\Yh_\bj(t_\ell)  \cdot 1/ n
 +
\sum_{\bi \in \I} \lambda_\bi^{1/2} \cdot
\lscp{B(\Xh(t_{\ell-1})) h_\bi}{h_\bj} \cdot
(\beta_\bi(t_\ell) - \beta_\bi (t_{\ell-1}))
\]
for $\bj \in \J$, where
\begin{equation}\label{g18}
\Xh (t) = \sum_{\bj \in \J} \Yh_\bj (t) \cdot h_\bj
\end{equation}
and
\begin{equation}\label{g19}
\Yh_\bj (0) = \scp{\xi}{h_\bj}.
\end{equation}
Equivalently, 
\[
\Yh_\bj (t_\ell) = \frac{1}{1+\mu_\bj / n }
\cdot \left( \Yh_\bj (t_{\ell-1})
+
\sum_{\bi \in \I} \lambda_\bi^{1/2} \cdot
\lscp{B(\Xh(t_{\ell-1})) h_\bi}{h_\bj} \cdot
(\beta_\bi(t_\ell) - \beta_\bi (t_{\ell-1}))\right).
\]

In general we define 
\[
0 = \tau_0 < \dots < \tau_M = 1
\]
by
\[
\{\tau_0, \dots, \tau_M\} = 
\bigcup_{\bi \in \I} \{ t_{0,\bi}, \dots, t_{n_\bi,\bi} \}.
\]
Moreover, we put
\[
\K_m = \{ \bi \in \I : \tau_m \in 
\{ t_{0,\bi}, \dots, t_{n_\bi,\bi} \} \}
\]
for $m=0,\dots,M$, and we define 
$s_{m,\bi}$ for $\bi \in \I$ and $m=1,\dots,M$ by
\[
s_{m,\bi} = \max ( 
\{ t_{0,\bi}, \dots, t_{n_\bi,\bi} \} \cap \left[0,\tau_m\right[).
\]
Finally, we use
\begin{equation}\label{semi}
\Gamma_\bj (t) = \prod_{\nu=1}^M
\frac{1}{1+\mu_\bj \cdot (t \wedge \tau_\nu - t \wedge
\tau_{\nu-1})}
\end{equation}
for approximation of the semigroup generated by $\Delta$. 
Then the drift-implicit Euler scheme is given by
\eqref{g18}, \eqref{g19}, and
\begin{multline}\label{f1}
\Yh_\bj (t) =  
\frac{\Gamma_\bj(t)}{\Gamma_\bj(\tau_{m-1})}\\
\mbox{} \cdot \left( \Yh_\bj (\tau_{m-1})
+
\sum_{\bi \in \K_m} \lambda_\bi^{1/2} \cdot
\lscp{B(\Xh(s_{m,\bi})) h_\bi}{h_\bj} \cdot
\frac{\Gamma_\bj (\tau_{m-1})}{\Gamma_\bj(s_{m,\bi})} \cdot
(\beta_\bi(\tau_m) - \beta_\bi (s_{m,\bi}))\right)
\end{multline}
for $\bj \in \J$, if
\[
t \in \left] \tau_{m-1},\tau_m\right].
\]
Equivalently, 
\begin{multline}\label{g20}
\Yh_\bj (t) = \Gamma_\bj (t) \cdot \scp{\xi}{h_\bj} \\
 +
\sum_{\bi \in \I} \lambda_\bi^{1/2}
\cdot \left( \sum_{t_{\ell,\bi} \leq \tau_m}
\lscp{B(\Xh(t_{\ell-1,\bi})) h_\bi}{h_\bj} \cdot
\frac{\Gamma_\bj(t)}{\Gamma_\bj(t_{\ell-1,\bi})} 
\cdot (\beta_\bi (t_{\ell,\bi}) - \beta_\bi (t_{\ell-1,\bi}))
\right). 
\end{multline}

For illustration we consider an example with
$\I = \{1,2\}$, $n_1 = 6$, and $n_2 = 4$.
Then, for instance, $\K_2 = \{2\}$, $\K_3=\{1\}$,
and $\K_4 = \{1,2\}$. Moreover, for $t \in \left]\tau_2,\tau_3\right]
= \left]1/4,1/3\right]$ the approximation $\Xh(t)$ is
based on the increments $\beta_1(1/6)$, $\beta_1(1/3) - \beta_1(1/6)$,
and $\beta_2(1/4)$, while $\beta_2(1/2) - \beta_2(1/4)$
is not used at all.

\section{Error Analysis}\label{s4}

Henceforth constants 
that are hidden in notations like%
\footnote{Suppose that $F$ and $G$ are functions on some set $A$
with values in $\left[0,\infty\right]$.
By definition, $F(a) \preceq G(a)$ means 
$F(a) \leq c \cdot G(a)$ for all $a \in A$ with some constant
$c \in \left]0,\infty\right[$. 
Furthermore, $F(a) \asymp G(a)$ means $F(a) \preceq G(a)$ and 
$G(a) \preceq F(a)$.}
$\preceq$ and $\asymp$ may only depend on $d$, $\xi$, $\lambda$ and $g$.

In the sequel we consider the particular choice 
\begin{equation}\label{g15}
\begin{aligned}
\I &= \{ \bi \in \N^d : |\bi|_2 \leq I\},\\
\J &= \{ \bj \in \N^d : |\bj|_2 \leq J\}
\end{aligned}
\end{equation}
in the definition of the approximations $\Xb$ and $\Xh$.
Then the error of the It\^o-Galerkin approximation
$\Xb$ is bounded as follows; see Section~\ref{sp} for the proof.

\begin{prop}\label{l4}
For $I, J > 0$
\[
E\left(\int_0^1  \|X(t) - \Xb(t)\|^2 \, dt\right) \preceq 
1/J^2 + \sum_{|\bi|_2 > I} \lambda_\bi / \mu_\bi.
\]
\end{prop}

Moreover, we have the following error bound for the implicit Euler 
scheme $\Xh$ with an arbitrary discretization \eqref{disc} specified by
$n \in \N^\I$; again we refer to Section \ref{sp} for the proof. 

\begin{theorem}\label{t1}
For $I, J > 0$ and $n \in \N^\I$
\[
E\left(\int_0^1  \|X(t) - \Xh(t)\|^2 \, dt \right)
\preceq 1/J^2 + 
\sum_{|\bi|_2 \leq I} \lambda_\bi / n_\bi + 
\sum_{|\bi|_2 > I} \lambda_\bi / \mu_\bi.
\]
\end{theorem}

Suppose that $\Xh$ may use a total of $N$ evaluations of
scalar Brownian motions $\beta_\bi$.
Then a proper choice of $I > 0$ and 
$n \in \N^\I$ is obtained
by minimizing 
\[
D(I,n) = \sum_{|\bi|_2 \leq I} \lambda_\bi / n_\bi + 
\sum_{|\bi|_2 > I} \lambda_\bi / \mu_\bi
\]
under the constraint $\sum_{|\bi|_2 \leq I} n_\bi \leq N$.
Up to a constant, the corresponding optimization problem is solved 
as follows.

Recall that, by assumption,
\begin{equation}\label{g46}
\lambda (r) = r^{-\gamma} \cdot L (r)
\end{equation}
with a slowly varying function $L : \left[1,\infty\right[ \to
\left]0,\infty\right[$. Let $N \in \N$.  We take 
\[
I = I_N = N^{1/(d+2)}
\]
to specify the scalar Brownian motions that
are evaluated by the algorithm.
For $\bi \in \N^d$ with $|\bi|_2 \leq I$ the Brownian motion
$\beta_\bi$ is evaluated with step-size $1/n_\bi$ where 
\[
n_\bi = n_{\bi,N} = \left\lceil \lambda_\bi^{1/2} \cdot
N^{ 1 - {\textstyle \frac{d-\gamma/2}{d+2}}} 
\cdot \left( L (N^{1/(d+2)}) \right)^{-1/2}\right\rceil
\]
if $\gamma \in \left[d,2d\right[$ and
\[
n_\bi = n_{\bi,N} = \left\lceil \lambda_\bi^{1/2} \cdot N \right\rceil
\]
if $\gamma \in \left]2d,\infty\right[$.
For the total number of evaluations we thus obtain
$\sum_{|\bi|_2 \leq I} n_\bi \asymp N$.
Moreover, 
\[
\inf \{ D(I,n) : I > 0,\ n \in \N^\I 
\ \text{with}\ \sum_{|\bi|_2 \leq I} n_\bi \leq N \}
\asymp D(I_N,n_N) \asymp e^2_*(N)
\]
with
\[
e_*(N) = N^{ -1/2 + {\textstyle \frac{d-\gamma/2}{d+2}}} 
\cdot \left( L (N^{1/(d+2)}) \right)^{1/2}
\]
if $\gamma \in \left[d,2d\right[$ and
\[
e_*(N) = N^{-1/2}
\]
if $\gamma \in \left]2d,\infty\right[$.
See M\"uller-Gronbach, Ritter (2006).

Finally, we take
\[
  J = J_N = e_*^{-1}(N).
\]
Hereby we have completely specified an algorithm
$\Xh = \Xh_N^*$.

\begin{theorem}\label{t2}
The error of the algorithm $\Xh_N^*$ satisfies
\[
\left(E\left( \int_0^1 \|X(t) - \Xh^*_N(t)\|^2 \, dt \right)\right)^{1/2}
\preceq e_*(N).
\]
\end{theorem}

\begin{rem}\label{2d}
The case of a regularly varying functions $\lambda$ of index 
$-2d$ is not covered by Theorem \ref{t2} 
but may be analyzed in a similar way.
Assume, for simplicity, that $\lambda(r) = r^{-2d}$.
Take $I_N$ as above, and define
\[
n_{\bi,N} = \lceil \lambda_\bi^{1/2} \cdot N/\ln N \rceil
\]
for $\bi\in\N^d$ with $|\bi|_2\le I_N$. Note that 
$\sum_{|\bi|_2 \leq I_N} n_{\bi,N} \asymp N$. Then
\[
\inf \{ D(I,n) : I > 0,\ n \in \N^\I 
\ \text{with}\ \sum_{|\bi|_2 \leq I} n_\bi \leq N \}
\asymp D(I_N,n_N) \asymp N^{-1}\cdot (\ln N)^2.
\]
Furthermore, take 
$J_N = N^{1/2} \cdot (\ln N)^{-1}$. Due to 
Theorem \ref{t1} the resulting algorithm $\Xh^*_N$ satisfies
\[
\left(E\left( \int_0^1 \|X(t) - \Xh^*_N(t)\|^2 \, dt\right) \right)^{1/2}
\preceq N^{-1/2}\cdot \ln N.
\] 
\end{rem}

\begin{rem}\label{uni1}
Consider the implicit Euler scheme $\Xh$ with a 
uniform time discretization
\eqref{disc}, i.e., $t_\ell = t_{\ell,\bi} = \ell/n$
for all $\bi\in \N^d$ with $|\bi|_2\le I$ and some constant $n\in\N$. 
Assume, for simplicity, that 
$\lambda(r) = r^{-\gamma}$ with $\gamma\in \left]d,\infty\right[$.
By Theorem \ref{t1},
\[
E\left( \int_0^1 \|X(t) - \Xh(t)\|^2 \, dt \right)\le 1/J^2 + d(I,n)
\]
with
\[
d(I,n) = 1/n \cdot \sum_{|\bi|_2\le I} |\bi|_2^{-\gamma} +
\sum_{|\bi|_2 > I} |\bi|_2^{-(\gamma + 2)}.
\]
Minimization of this quantity, up to a constant, 
subject to the constraint $n\cdot \# \I \le N$ leads to
\[
\inf \{ d(I,n) : I > 0,\ 
n \in \N \ \text{with}\ n\cdot \# \I\leq N \}
\asymp d(I_N,n_N) \asymp N^{- 1 + {\textstyle \frac{d}{\gamma+2}}}
\]
with
\[
I = I_N = N^{1/(\gamma+2)}
\]
and
\[
n = n_N = 
\left\lceil N^{(\gamma+2-d)/(\gamma+2)}\right\rceil.
\]
Take 
\[
J = J_N = N^{1/2 - {\textstyle \frac{d}{2(\gamma+2)}}},
\]
and let $\Xh^{\text{uni}}_N$ denote the resulting algorithm.
By definition, 
$n_N\cdot \# \I_N \asymp N$ for the total
number of evaluations of scalar Brownian motions $\beta_\bi$ used by 
$\Xh^{\text{uni}}_N$, and Theorem \ref{t1} yields
\begin{equation}\label{n1}
\left( E \left(\int_0^1 \|X(t) - \Xh^{\text{uni}}_N(t)\|^2 \, dt 
\right)\right)^{1/2}
\preceq N^{-1/2 + {\textstyle \frac{d}{2(\gamma+2)}}}.
\end{equation}
\end{rem}

\begin{rem}\label{vgl}
We compare the implicit Euler schemes $\Xh^*_N$
and $\Xh^{\text{uni}}_N$, both of which 
roughly use $N$ evaluations of scalar Brownian motions.
Assume that $\lambda(r) = r^{-\gamma}$ with 
$\gamma\in \left]d,\infty\right[\setminus \{2d\}$, and put
\[
\alpha^*(\gamma,d) = 
\frac{1}{2}-\frac{(2d-\gamma)_+}{2(d+2)},
\qquad \alpha(\gamma,d) = \frac{1}{2}-\frac{d}{2(\gamma+2)}.
\]
{}From Theorem \ref{t2} and Remark \ref{uni1} we get
\[
\left(E \left( \int_0^1 \|X(t) - \Xh^*_N(t)\|^2 \, dt \right)\right)^{1/2}
\preceq N^{-\alpha^*(\gamma,d)}
\]
for the non-uniform discretization, and
\[
\left( E \left(\int_0^1 \|X(t) - \Xh^{\text{uni}}_N(t)\|^2 \, dt 
\right)\right)^{1/2}
\preceq N^{-\alpha(\gamma,d)}
\]
for the uniform discretization. 

We always have 
\[
\alpha^*(\gamma,d) > \alpha(\gamma,d).
\]
In the limit for a low degree of smoothness 
\[
\lim_{\gamma\to d+} \alpha^*(\gamma,d) = 
\lim_{\gamma\to d+} \alpha(\gamma,d) = 1/(d+2).
\]
Conversely, for a high degree of smoothness
\[
\lim_{\gamma\to \infty} \alpha(\gamma,d) = 1/2,
\]
while $\alpha^*(\gamma,d) =1/2$ already holds if $\gamma > 2d$.

\end{rem}

\section{Optimality}\label{s5}

The results from Section \ref{s4} provide upper bounds
for the error of specific algorithms. 
In particular, the comparison of the implicit Euler schemes
based on uniform and non-uniform time discretizations
is in fact a comparison only of the corresponding upper bounds. 
It is therefore important to know whether these
upper bounds are lower bounds for the error as well, and, even more, 
to raise the following questions:
\begin{itemize}
\item[(i)]
Does there exist any algorithm $\Xh_N$ that
uses a total of $N$ evaluations of
scalar Brownian motions $\beta_\bi$ and achieves an error 
significantly smaller than the upper bound $e_*(N)$ for the 
algorithm $\Xh^*_N$?
\item[(ii)] 
Are non-uniform time discretizations superior to uniform ones?
\end{itemize}

To answer these questions we consider arbitrary methods that evaluate a 
finite number of 
Brownian motions $\beta_\bi$ at a finite number of points and then 
produce a curve in $H$ that is close to the corresponding realization 
of $X$.
In general, the selection and evaluation of the scalar Brownian motions
$\beta_\bi$, is specified by a finite set
\[
\I \subseteq \N^d
\]
and nodes
\[
0 < t_{1,\bi} < \dots < t_{n_\bi,\bi} \leq 1
\]
for $\bi \in \I$ and $n_\bi \in \N$. 
Every Brownian motion $\beta_\bi$ with $\bi \in \I$ is evaluated
at the corresponding nodes $t_{\ell,\bi}$.
The total number of evaluations is given by
\[
|n|_1 = \sum_{\bi\in\I} n_\bi.
\]
Formally, an approximation $\Xh$ to $X$ is given by
\begin{equation}\label{yu}
\Xh(t) = \phi \bigl(t, 
\beta_{\bi_1}(t_{1,\bi_1}), \dots, \beta_{\bi_1}(t_{n_{\bi_1},\bi_1}),
\dots  ,
\beta_{\bi_k}(t_{1,\bi_k}), \dots,
\beta_{\bi_k}(t_{n_{\bi_k},\bi_k})\bigr),
\end{equation}
where 
\[
\phi : [0,1] \times \R^{|n|_1} \to H
\]
is any measurable mapping and $\I = \{\bi_1, \dots , \bi_k\}$.
Here $\phi$ may depend in any way on the initial value 
$\xi$, the eigenvalues $\lambda_\bi$, and the function $g$,
which is used to define the mapping $B$ in the heat equation
\eqref{g1}. The error of $\Xh$ is defined by
\[
e(\Xh) = \left(E \left(
 \int_0^1 \| X(t) - \Xh(t) \|^2 \, dt \right) \right)^{1/2},
\]
cf.\ Theorems \ref{t1} and \ref{t2} and the subsequent Remarks.

Let $\XX_N$ denote the class of all algorithms \eqref{yu}
that use a total of $N$
evaluations of the
scalar Brownian motions $\beta_\bi$, i.e., $|n|_1 = N$. 
We wish to minimize the error in this class,  
and hence we study the $N$th minimal error
\[
e(N) = \inf_{\Xh \in \XX_N} e(\Xh).
\]

In particular, our algorithm $\Xh^*_N$ is of the form \eqref{yu},
and its total number of evaluations of scalar Brownian motions
is roughly given by $N$.

We obtain a negative answer to Question (i). 

\begin{theorem}\label{t3}
The sequence of algorithms $\Xh^*_N$
is asymptotically optimal, i.e., 
\[
e(\Xh_N^*) \asymp e(N),
\]
and
\[
e(N) \asymp e_*(N).
\]
\end{theorem}

\begin{proof}
In view of Theorem \ref{t2} it remains to show that
$e(N) \succeq e_*(N)$, and this lower bound is
is a consequence of a more general result 
established in M\"uller-Gronbach, Ritter (2006, Thm.\ 1).
\end{proof}

With respect to Question (ii) one needs to study the subclass 
$\XX_N^{\text{uni}}\subset \XX_N$ of algorithms that are based on 
a uniform discretization, i.e., 
$t_{\ell,\bi} = \ell / n$ for all $\bi\in\I$ and 
some constant $n\in\N$.
The corresponding $N$th minimal error in this class is given by
\[
e^{\text{uni}}(N) = \inf_{\Xh \in \XX_N^{\text{uni}}} e(\Xh).
\]  

\begin{rem}
Consider the specific equation
\begin{equation}\label{n2}
\begin{aligned}
dX(t) &= \Delta X(t) \, dt +  dW(t),\\
X(0) &= 0,
\end{aligned}
\end{equation}
i.e., $g=1$ or, equivalently, $B(x)= \text{id}$, and assume that
$\lambda(r)=r^{-\gamma}$ with 
$\gamma\in\left] d,\infty\right[\setminus \{2d\}$. Then
\begin{equation}\label{n3}
e^{\text{uni}}(N) \succeq N^{-1/2+{\textstyle \frac{d}{2(\gamma+2)}}}
\end{equation}
see M\"uller-Gronbach, Ritter (2006, Remark 6), so that
\[
e(\Xh^{\text{uni}}_N) \asymp 
e^{\text{uni}}(N) \asymp N^{-1/2+{\textstyle \frac{d}{2(\gamma+2)}}}
\]
follows from Remark \ref{uni1}.

We thus conclude that the upper bound \eqref{n1} for the error
of the implicit Euler scheme $\Xh^{\text{uni}}_N$ is sharp.
Moreover, these algorithms form an asymptotically optimal
sequence among all algorithms that use uniform time discretizations.
Our comparison of orders of convergence in Remark \ref{vgl}
is therefore a result on minimal errors and clearly shows
the superiority of non-uniform time discretizations.
\end{rem}

Note that the conclusions from the previous remark only apply
to the specific equation \eqref{n2}. We conjecture, however,
that the lower bound \eqref{n3} 
holds in general, in which case these conclusions hold in general
as well.

\section{Proofs}\label{sp}

We start with regularity properties of the solution $X$ of 
equation \eqref{g1}. 
For the mean-square smoothness of $X$ we have
\begin{equation}\label{g4}
E\| X(s) - X(t)\|^2 \preceq |t-s| \cdot ( 1 + \psi(\min(s,t)))
\end{equation}
for $s,t \in [0,1]$ with
\[
\psi(t) = \sum_{\bi \in \N^d} \mu_\bi \cdot E(\scp{X(t)}{h_\bi}^2)
\]
satisfying
\begin{equation}\label{g10}
\int_0^1 \psi(t)\, dt < \infty. 
\end{equation}
See M\"uller-Gronbach, Ritter (2006, Lemma 1).

Consider the Sobolev space $W^1_2 = W^1_2(\left]0,1\right[^d)$
and its subspace $W^{1,0}_2$.
Note that $h_\bi \in W^{1,0}_2$ and
\[
\scp{x}{h_\bi}_{W^1_2} = (1+\mu_\bi) \cdot \scp{x}{h_\bi}
\]
for every $x \in W^{1,0}_2$. Consequently, the functions
$(1+\mu_\bi)^{-1/2} \cdot h_\bi$ form a 
complete orthonormal system in $W^{1,0}_2$, and
\begin{equation}\label{g9}
W^{1,0}_2 = \{ x \in H : 
\sum_{\bi \in \N^d} \mu_\bi \cdot \scp{x}{h_\bi}^2 < \infty\}
\end{equation}
as well as
\begin{equation}\label{g8}
\|x\|^2_{W^1_2} \leq 2 \cdot
\sum_{\bi \in \N^d} \mu_\bi \cdot \scp{x}{h_\bi}^2
\end{equation}
for $x \in W^{1,0}_2$, which is the Poincar\'e inequality.

\begin{lemma}\label{l3}
For Lebesgue-almost every $t\in[0,1]$ we have
\[
X(t) \in W^{1,0}_2 
\]
with probability one and
\[
\sum_{\bj\in\N^d} 1/\mu_\bj \cdot 
E\scp{T_g(X(t))\cdot h_\bi}{h_\bj}^2
\preceq 1/\mu_\bi \cdot (1+E\|X(t)\|_{W^1_2}^2)
\]
for every $\bi\in\N^d$. Moreover,
\[
\int_0^1 E \| X(t)\|^2_{W^{1}_2} \, dt < \infty.
\]
\end{lemma}

\begin{proof}
Combine \eqref{g10}, \eqref{g9}, and \eqref{g8} to obtain the
first and the last claim.

For the proof of the second claim we
note that 
\begin{equation}\label{g7}
\|T_g (x)\|_{W^1_2} \preceq 1 +  \|x\|_{W^1_2}
\end{equation}
for $x \in W^1_2$, see Appell, Zabrejko (1990, Theorems 9.2 and 9.5).
Furthermore, we may assume $g(0)=0$ without loss of generality. Then
\begin{equation}\label{g6}
T_g(W^{1,0}_2) \subset W^{1,0}_2
\end{equation}
is easily verified.
In view of \eqref{g7}, \eqref{g6} and the first statement in the lemma
it suffices to show that

\begin{equation}\label{h1}
\sum_{\bj \in \N^d} 1/\mu_\bj \cdot \scp{x\cdot h_\bi}{h_\bj}^2
\preceq 1/\mu_\bi \cdot \|x\|_{W^1_2}^2
\end{equation}
for all $x \in W^{1,0}_2$ and $\bi \in \N^d$.

To this end fix $\bi, \bj \in \N^d$ and $\ell \in \{1,\dots,d\}$, 
and put
\[
f_\bi = 1/(i_\ell \pi) \cdot \tfrac{\partial}{\partial u_\ell} h_\bi,
\qquad
f_\bj = 1/(j_\ell \pi) \cdot \tfrac{\partial}{\partial u_\ell} h_\bj.
\]
Then
\begin{align*}
i_\ell^2 \cdot \scp{x}{h_\bi \cdot h_\bj}^2 
& \asymp 
\scp{x}{\tfrac{\partial}{\partial u_\ell} f_\bi \cdot h_\bj}^2\\
&=
\left(
\scp{\tfrac{\partial}{\partial u_\ell}x}{f_\bi \cdot h_\bj} -
j_\ell \pi \cdot \scp{x}{f_\bi \cdot f_\bj}  \right)^2 \\
&\preceq
\scp{\tfrac{\partial}{\partial u_\ell}x}{f_\bi \cdot h_\bj}^2 +
j_\ell^2 \cdot \scp{x}{f_\bi \cdot f_\bj}^2.
\end{align*}
Hereby
\[
i_\ell^2 \cdot \sum_{\bj \in \N^d} 1/\mu_\bj \cdot \scp{x}{h_\bi\cdot
h_\bj}^2
\preceq 
\left\|\tfrac{\partial}{\partial u_\ell}x \cdot f_\bi \right\|^2 +
\left\|x \cdot f_\bi \right\|^2,
\]
and we conclude that
\[
\mu_\bi \cdot \sum_{\bj \in \N^d} 1/\mu_\bj \cdot 
\scp{x\cdot h_\bi}{h_\bj}^2 \preceq \|x\|^2_{W^1_2},
\]
which yields \eqref{h1}.
\end{proof}

\subsection{Properties of the It\^{o}-Galerkin approximation}

Let $P_\I$ and $P_\J$ denote the orthogonal projections
onto the subspaces $\spann \{h_\bi: \bi \in \I\}$ and
$\spann \{h_\bj: \bj \in \J\}$, respectively, and put
\[
\Bb (x) = P_\J \circ B(x) \circ P_\I.
\]
Then $\Bb: H\to \LL$ satisfies \eqref{e1} and	
$\Xb$ is the mild solution of \eqref{g1} with $B$ being
replaced by $\Bb$. Hence
\begin{equation}\label{g14}
\sup_{t \in [0,1]} E\| \Xb(t)\|^2 \leq c_1, 
\end{equation}
see \eqref{g3}.

We establish an error bound for piecewise constant interpolation of
$\Xb$. 

\begin{lemma}\label{l5}
For $\I, \J \subset \N^d$ and $m \in \N$
\[
\sum_{\ell=0}^{m-1} 
\int_{\ell/m}^{(\ell+1)/m} E \|\Xb(t) - \Xb(\ell/m)\|^2 \, dt 
\preceq 1/m.
\]
\end{lemma}

\begin{proof}
Note that \eqref{g4} and \eqref{g10} are valid, too, for
$\Xb$ and
\[
\pb (t) = \sum_{\bj\in\J} \mu_\bj \cdot E(\Yb_\bj^2(t))
\]
instead of $X$ and $\psi$, respectively.
For $\ell \in \{0,\dots,m-1\}$
take $s_\ell \in [\ell/m,(\ell+1)/m]$ with
\[
\pb (s_\ell) / m \leq  \int_{\ell/m}^{(\ell+1)/m} \pb(t)\, dt.
\]

On the first subinterval,
\[
\int_{0}^{1/m} E \|\Xb(t) - \Xb(0)\|^2 \, dt \leq
2/m \cdot \sup_{t \in [0,1]} E\| \Xb(t)\|^2 \preceq 1/m,
\]
see \eqref{g14}.
On the subintervals $[\ell/m,(\ell+1)/m]$ with $\ell \ge 1$ we
proceed as follows.
If $t \in [\ell/m,s_\ell]$, then
\begin{align*}
 E \|\Xb(t) - \Xb(\ell/m)\|^2 
& \preceq
E \|\Xb(t) - \Xb(s_{\ell-1})\|^2  +
E \|\Xb(s_{\ell-1}) - \Xb(\ell/m)\|^2\\
&\preceq 1/m \cdot (1 + \pb(s_{\ell-1}) )\\
&\preceq 1/m + \int_{(\ell-1)/m}^{\ell/m} \pb (s) \, ds.
\end{align*}
If $t \in [s_\ell,(\ell+1)/m]$, then
\begin{align*}
& E \|\Xb(t) - \Xb(\ell/m)\|^2 \\
& \qquad \preceq
E \|\Xb(t) - \Xb(s_\ell)\|^2  +
E \|\Xb(s_\ell) - \Xb(s_{\ell-1})\|^2  +
E \|\Xb(s_{\ell-1}) - \Xb(\ell/m)\|^2\\
&\qquad\preceq 1/m \cdot (1 + \pb (s_\ell) + \pb(s_{\ell-1}) )\\
&\qquad\preceq 1/m + \int_{(\ell-1)/m}^{(\ell+1)/m} \pb (s) \, ds.
\end{align*}
We conclude that
\[
\int_{\ell/m}^{(\ell+1)/m} E \|\Xb(t) - \Xb(\ell/m)\|^2 \, dt 
\preceq 1/m^2 + 1/m \cdot 
\int_{(\ell-1)/m}^{(\ell+1)/m} \pb (s) \, ds,
\]
which completes the proof.
\end{proof}

\begin{proof}[Proof of Proposition \ref{l4}]
Recall the particular choice \eqref{g15} of the sets $\I$ and $\J$
and let $c$ denote the right-hand side in Proposition \ref{l4}.
Moreover, let
\[
X^{(k)}(t) = \sum_{|\bj|_2 \leq J} Y^{(k)}_\bj(t) \cdot h_\bj
\]
for $k=1,2$ with
\[
Y^{(1)}_\bj(t) =
\sum_{|\bi|_2 > I} \lambda_\bi^{1/2} \cdot 
\int_0^t \exp(-\mu_\bj (t-s)) \cdot \lscp{T_g (X(s)) \cdot h_\bi}
{h_\bj} \, d\beta_\bi(s)
\]
and  
\[
Y^{(2)}_\bj(t) = \exp(- \mu_\bj t) \cdot \scp{\xi}{h_\bj} +
\sum_{|\bi|_2 \leq I} \lambda_\bi^{1/2} \cdot 
\int_0^t \exp(-\mu_\bj (t-s)) \cdot \lscp{T_g (X(s)) \cdot h_\bi}
{h_\bj} \, d\beta_\bi(s).
\]
Then
\[
X(t) = X^{(1)}(t) + X^{(2)}(t) + \sum_{|\bj|_2 > J} Y_\bj(t)\cdot h_\bj,
\]
and consequently
\begin{multline*}
\int_0^t E \|X(s) - \Xb(s)\|^2 \, ds \\
\preceq 
\sum_{|\bj|_2 > J} \int_0^1 E(Y_\bj^2(t)) \, dt
+
\int_0^1 E \|X^{(1)} (t) \|^2 \, dt 
+
\int_0^t E \|X^{(2)}(s) - \Xb(s)\|^2 \, ds.
\end{multline*}

We have
\[
\int_0^1 E\bigl(Y^{(1)}_\bj(t)\bigr)^2 \, dt
\preceq \sum_{|\bi|_2 > I} \lambda_\bi / \mu_\bj
\cdot \int_0^1 E\lscp{T_g (X(t)) \cdot h_\bi}{h_\bj}^2 \, dt,
\] 
and therefore
\[
\int_0^1 E \|X^{(1)} (t) \|^2 \, dt \leq
\sum_{\bj \in \N^d} 
\int_0^1 E(Y^{(1)}_\bj(t))^2 \, dt
\preceq 
\sum_{|\bi|_2 > I} \lambda_\bi / \mu_\bi \leq c
\]
by Lemma \ref{l3}. Furthermore, 
\[
\sum_{|\bj|_2 \geq J} \int_0^1 E(Y_\bj^2(t)) \, dt \preceq 1/J^2 \leq
c
\]
follows from \eqref{g13}, \eqref{g3},
and $\sup_{\bi \in \N^d} \| h_i\|_\infty < \infty$.
Finally, if $|\bj|_2 \leq J$, then
\[
E \bigl(Y^{(2)}_\bj(t) - \Yb_\bj(t)\bigr)^2 
\leq
\sum_{|\bi|_2 \leq I} \lambda_\bi \cdot
\int_0^t 
E \lscp{T_g (X(s)) - T_g (\Xb(s))}{h_\bi\, h_\bj}^2  \, ds,
\]
and due to \eqref{g13} 
\begin{align*}
E \|X^{(2)}(t) - \Xb(t)\|^2  &\preceq
\int_0^t  E \|T_g (X(s)) - T_g (\Xb(s))\|^2  \, ds
\preceq
\int_0^t  E \|X(s) - \Xb(s)\|^2 \, ds\\
&\preceq 2c + 
\int_0^t E \|X^{(2)}(s) - \Xb(s)\|^2\, ds.
\end{align*}
Since $E \lscp{X^{(1)}(t)}{X^{(2)}(t)} = 0$, we get
$E\|X^{(2)}(t)\|^2 \leq E\|X(t)\|^2$. Use \eqref{g3} and \eqref{g14} 
to conclude that
\[
\sup_{t \in [0,1]}  
E \|X^{(2)}(t) - \Xb(t)\|^2 < \infty.
\]
It remains to apply Gronwall's Lemma to complete the proof.
\end{proof}

\subsection{Properties of the implicit Euler scheme}

Recall the definition \eqref{semi} of $\Gamma_\bj$ used for
approximation of the semigroup.

\begin{lemma}\label{l7}
Suppose that $\bi\in \I$  and $\bj\in \J$.
Then, for $\ell = 0,\dots,n_\bi-1$,
\[
\int_{t_{\ell,\bi}}^1
\frac{\Gamma_\bj^2(t)}{\Gamma_\bj^2(t_{\ell,\bi})} \, dt 
\leq 2/\mu_\bj
\]
as well as
\[
\int_{t_{\ell,\bi}}^1
\left( \frac{\Gamma_\bj(t)}{\Gamma_\bj(t_{\ell,\bi})} -
\exp(-\mu_\bi (t-t_{\ell,\bi})) \right)^2 \, dt 
\preceq 1/n^*,
\]
where
\[
n^* = \max \{ n_\bi: \bi\in\I\}.
\]
Furthermore, for $0 \leq s \leq t \leq 1$, 
\[
\left|1- \frac{\Gamma_\bj(t)}{\Gamma_\bj(s)}
\right| \leq \min(1, \mu_\bj \cdot (t-s)).
\]
\end{lemma}

\begin{proof}
For $t \in [t_{k,\bi},t_{k+1,\bi}]$ with $k \geq \ell$
\[
\frac{\Gamma_\bj(t)}{\Gamma_\bj(t_{\ell,\bi})}
\leq \frac{1}{(1+\mu_\bj/n_\bi)^{k-\ell}} \cdot 
\frac{1}{1+\mu_\bj \cdot (t - t_{k,\bi})},
\]
and therefore
\[
\int_{t_{\ell,\bi}}^1
\frac{\Gamma_\bj^2(t)}{\Gamma_\bj^2(t_{\ell,\bi})} \, dt 
\leq \frac{1}{\mu_\bj + n_\bi} \cdot
\sum_{k=0}^{n_\bi-1} 
\frac{1}{(1+\mu_\bj/n_\bi)^{2k}}.
\]
Thus, if $\mu_\bj / n_\bi \geq 1$,
\[
\int_{t_{\ell,\bi}}^1
\frac{\Gamma_\bj^2(t)}{\Gamma_\bj^2(t_{\ell,\bi})} \, dt 
\leq 2/\mu_\bj,
\]
and otherwise
\[
\int_{t_{\ell,\bi}}^1
\frac{\Gamma_\bj^2(t)}{\Gamma_\bj^2(t_{\ell,\bi})} \, dt 
\leq \frac{1}{n_\bi} \cdot
\frac{1}{1- 1/(1+ \mu_\bj / n_\bi)^2} \leq 2/\mu_\bj,
\]
too.

For the proof of the second statement put 
$k^* = \lceil t_{\ell,\bi} \cdot n^* \rceil$ and
\[
f(t) = \frac{\Gamma_\bj(t)}{\Gamma_\bj(t_{\ell,\bi})} -
\exp(-\mu_\bj (t-t_{\ell,\bi})). 
\]
Then $0 \leq f \leq 1$ and
\[
\int_{t_{\ell,\bi}}^1
f^2(t) \, dt 
\preceq 1/n^* + \sum_{k=k^*}^{n^*-1}
\int_{k/n^*}^{(k+1)/n^*} f^2(t) \, dt. 
\]
It remains to show that
\begin{equation}\label{g17}
\sum_{k=k^*}^{n^*-1} \sup_{t \in [k/n^*,(k+1)/n^*]} f^2(t) \preceq 1.
\end{equation}

To this end 
assume that $t \in [k/n^*,(k+1)/n^*]$ for some $k \geq k^*$ in the 
sequel. Use
\[
\frac{\Gamma_\bj(t)}{\Gamma_\bj(t_{\ell,\bi})}
\leq
\frac{1}{1+\mu_\bj \cdot (k^*/n^* - t_{\ell,\bi})} \cdot
\frac{1}{(1+\mu_\bj/n^*)^{k-k^*}} \cdot
\frac{1}{1+\mu_\bj \cdot (t - k/n^*)}
\]
to obtain
\[
f(t) 
\leq
\left(\frac{1}{(1+\mu_\bj/n^*)^{k-k^*}} - \exp(-\mu_\bj (k-k^*)/n^*)\right)
+
\exp(-\mu_\bj (k-k^*)/n^*) \cdot 
\left( f_0 + f_1(t)
\right)
\]
with
\[
f_0 =
\frac{1}{1+\mu_\bj \cdot (k^*/n^* - t_{\ell,\bi})} -
\exp(-\mu_\bj (k^*/n^* - t_{\ell,\bi})) 
\]
and
\[
f_1(t) =
\frac{1}{1+\mu_\bj \cdot (t-k/n^*)} -
\exp(-\mu_\bj (t-k/n^*)). 
\]
Note that $1/(1+u) - \exp(-u) \preceq \min(1/(1+u),u^2)$ for $u \geq 0$.
Let $u = \mu_\bj / n^*$.
Then
\[
\frac{1}{(1+\mu_\bj/n^*)^{k-k^*}} - \exp(-\mu_\bj (k-k^*)/n^*)
\preceq
\frac{1}{(1+u)^{k-k^*-1}} \cdot \min(1/(1+u),u^2)
\]
for $k > k^*$, and therefore
\[
\sum_{k=k^*}^{n^*-1} 
\left(\frac{1}{(1+u)^{k-k^*}} - \exp(- (k-k^*) \cdot u )\right)^2
\preceq
\frac{(1+u)^2}{(1+u)^2 - 1} \cdot \min(1/(1+u)^2,u^4)
\preceq 1.
\]
Since $\max(f_0,f_1(t)) \preceq \min (1, u^2)$, we have
\[
\sum_{k=k^*}^{n^*-1} 
\exp(-2\mu_\bj (k-k^*)/n^*) \cdot (f_0 + f_1(t))
\preceq 
\frac{1}{ 1 - \exp (-2u)} \cdot \min(1,u^2) \preceq 1,
\]
which completes the proof of \eqref{g17}.

For the proof of the third statement let $s \leq t$,
and assume that
$s \in [\tau_{\kappa-1},\tau_\kappa]$ and 
$t \in [\tau_{\nu-1},\tau_\nu]$. By definition
\[
\frac{\Gamma_\bj(t)}{\Gamma_\bj(s)} 
= 
\frac{1 + \mu_\bj \cdot (s-\tau_{\kappa-1})}
     {1 + \mu_\bj \cdot (\tau_\kappa-\tau_{\kappa-1})}
\cdot 
\prod_{\iota = \kappa+1}^{\nu-1} 
\frac{1}{1+\mu_\bj \cdot (\tau_\iota - \tau_{\iota-1})}
\cdot \frac{1}{1+\mu_\bj \cdot (t - \tau_{\nu-1})},
\]
which implies
\begin{align*}
\left| 1 - \frac{\Gamma_\bj(t)}{\Gamma_\bj(s)} \right|
& \leq
\left| 1 - 
\frac{1 + \mu_\bj \cdot (s-\tau_{\kappa-1})}
     {1 + \mu_\bj \cdot (\tau_\kappa-\tau_{\kappa-1})} \right| \\
&
\phantom{=} \mbox{} + 
\sum_{\iota = \kappa+1}^{\nu-1} 
\left| 1 - \frac{1}{1+\mu_\bj \cdot (\tau_\iota - \tau_{\iota-1})}
\right|
+
\left| 1 - \frac{1}{1+\mu_\bj \cdot (t - \tau_{\nu-1})}\right|\\
&\leq 
\mu_\bj \cdot (t - s).
\end{align*}
Finally, $0 < \Gamma_\bj(t)/\Gamma_\bj(s) \leq 1$.
\end{proof}

Put
\[
\at_{\bi,\bj}(t) = 
E\left(\lscp{T_g(\Xh(t)) \cdot h_\bi}{h_\bj}^2\right),
\]
and note that
\begin{equation}\label{g26}
\sum_{\bj \in \J} \at_{\bi,\bj}(t) \preceq 1 + E \|\Xh(t)\|^2
\end{equation}
due to \eqref{g13} and $\sup_{\bi \in \N^d} \|h_\bi\|_\infty \preceq 1$.

\begin{lemma}\label{l6}
For $\I, \J \subset \N^d$ and $n \in \N^{\I}$ 
\[
\sup_{t \in [0,1]} E\| \Xh(t)\|^2 \preceq 1. 
\]
\end{lemma}

\begin{proof}
At first we slightly modify the process $\Xh$
by replacing the Brownian increments 
$\beta_\bi(\tau_m) - \beta_\bi (s_{m,\bi})$
in the definition \eqref{f1} of $\Yh_\bj$ by
increments $\beta_\bi(t) - \beta_\bi (s_{m,\bi})$.
More precisely, we consider 
$\Xt (t) = \sum_{\bj \in \J} \Yt_\bj (t) \cdot
h_\bj$ with $\Yt_\bj (0) = \scp{\xi}{h_\bj}$ and
\begin{multline*}
\Yt_\bj (t) =  
\frac{\Gamma_\bj(t)}{\Gamma_\bj(\tau_{m-1})} \\
\mbox{} \cdot \left( \Yt_\bj (\tau_{m-1})
+
\sum_{\bi \in \K_m} \lambda_\bi^{1/2} \cdot
\lscp{B(\Xt(s_{m,\bi})) h_\bi}{h_\bj} \cdot
\frac{\Gamma_\bj (\tau_{m-1})}{\Gamma_\bj(s_{m,\bi})} \cdot
(\beta_\bi(t) - \beta_\bi (s_{m,\bi}))\right)
\end{multline*}
for $t \in \left]\tau_{m-1},\tau_m\right]$.
Note that $\Yh_\bj$ and $\Yt_\bj$ as well as 
$\Xh$ and $\Xt$ coincide at the points $\tau_m$.
Moreover, by construction of these processes we have
\begin{equation}\label{hilfe}
\Yt_\bj(\tau_m) \text{ and } \Xt(\tau_m) \text{ are measurable w.r.t. }
\sigma\left(\{\beta_\bi(t_{\ell,\bi}):\ 
t_{\ell,\bi}\le \tau_m, \bi\in\I\}\right).
\end{equation}
We claim that
\begin{equation}\label{g24}
\sup_{t \in [0,1]} E\| \Xt(t)\|^2 \preceq 1. 
\end{equation}

Assume that
$t \in \left] \tau_{m-1},\tau_m\right]$ in the following.
Observing \eqref{hilfe} we obtain
\begin{align*}
&E (\Yt_\bj(t) - \Yt_\bj(\tau_{m-1}))^2 \\
& =  
\left( 1 - \frac{\Gamma_\bj(t)}{\Gamma_\bj(\tau_{m-1})} \right)^2 \cdot
E (\Yt_\bj^2(\tau_{m-1})) 
+
\sum_{\bi \in \K_m} \lambda_\bi \cdot
\at_{\bi,\bj} (s_{m,\bi}) \cdot
\frac{\Gamma^2_\bj (\tau_{m-1})}{\Gamma^2_\bj(s_{m,\bi})} \cdot
(t - s_{m,\bi})\\
&\leq 
E (\Yt_\bj^2(\tau_{m-1})) + 
\sum_{\bi \in \K_m} \lambda_\bi/n_\bi \cdot
\at_{\bi,\bj} (s_{m,\bi}). 
\end{align*}
{}From \eqref{g26} we therefore get
\[
E\|\Xt(t) -\Xt(\tau_{m-1})\|^2
\preceq  
1 + \max_{k=0,\dots,m-1} E \|\Xt(\tau_k)\|^2,
\]
and we conclude that
\[
f(s) = \sup_{r \in [0,s]} E \|\Xt(r)\|^2
\]
is finite for $s \in  [0,1]$, since $E\|\Xt(0)\|^2= \|\xi\|^2 <
\infty$.

Analogously to \eqref{g20} we have
\begin{multline*}
\Yt_\bj (t) = \Gamma_\bj (t) \cdot \scp{\xi}{h_\bj} \\
 +
\sum_{\bi \in \I} \lambda_\bi^{1/2}
\cdot \left( \sum_{t_{\ell,\bi} \leq \tau_m}
\lscp{B(\Xt(t_{\ell-1,\bi})) h_\bi}{h_\bj} \cdot
\frac{\Gamma_\bj(t)}{\Gamma_\bj(t_{\ell-1,\bi})} 
\cdot (\beta_\bi (t \wedge t_{\ell,\bi}) - \beta_\bi (t_{\ell-1,\bi}))
\right),
\end{multline*}
which implies
\[
E(\Yt_\bj^2(t))
=
\Gamma_\bj^2(t) \cdot \scp{\xi}{h_\bj}^2 
+
\sum_{\bi \in \I} \lambda_\bi \cdot
\left( \sum_{t_{\ell,\bi} \leq \tau_m} \at_{\bi,\bj}( t_{\ell-1,\bi})
\cdot \frac{\Gamma_\bj^2(t)}{\Gamma_\bj^2(t_{\ell-1,\bi})} \cdot
(t \wedge t_{\ell,\bi} - t_{\ell-1,\bi})
\right)
\]
due to the measurability property \eqref{hilfe}.
Use \eqref{g26} to derive
\begin{align*}
E\|\Xt (t)\|^2
&\preceq
\|\xi\|^2 + 
\sum_{\bi in \I} \lambda_\bi
\cdot \left( \sum_{t_{\ell,\bi} \leq \tau_m} (1 + f(t_{\ell-1,\bi}))
\cdot
(t \wedge t_{\ell,\bi} - t_{\ell-1,\bi})
\right) \\
&\preceq
1 + \int_0^t f(s) \, ds,
\end{align*}
so that \eqref{g24} follows by means of Gronwall's Lemma.

For the process $\Xh$ we apply
\eqref{g20} and observe \eqref{hilfe} again to obtain 
\begin{equation}\label{g23}
E(\Yh_\bj^2(t))
=
\Gamma_\bj^2(t) \cdot \scp{\xi}{h_\bj}^2 
+
\sum_{\bi \in \I} \lambda_\bi/n_i \cdot
\left( \sum_{t_{\ell,\bi} \leq \tau_m} 
\ah_{\bi,\bj}(t_{\ell-1,\bi})
\cdot \frac{\Gamma_\bj^2(t)}{\Gamma_\bj^2(t_{\ell-1,\bi})} 
\right).
\end{equation}
Using \eqref{g24} we conclude that
\[
E\|\Xh (t)\|^2
\preceq
\|\xi\|^2 + 
\sum_{\bi \in \I} \lambda_\bi
\cdot 
\left( 
1 + \max_{\ell=0,\dots,n_\bi} E\|\Xt(t_{\ell,\bi})\|^2
\right)
\preceq
1. 
\]
\end{proof}

We turn to the mean-square regularity of the process $\Xh$.

\begin{lemma}\label{l8}
For $\I, \J \subset \N^d$, $n \in \N^{\I}$, and
$0 \leq s \leq t \leq 1$
\[
E\| \Xh(s) - \Xh(t)\|^2 
\preceq
(t-s) \cdot (1 + \ph (s)) + \sum_{\bi \in \I} \lambda_\bi / n_\bi,
\]
where
\[
\ph (s) = \sum_{\bj \in \J} \mu_\bj \cdot E(\Yh^2_\bj(s)).
\]
Moreover,
\begin{equation}\label{g22}
\int_0^1 \ph(s) \, ds \preceq 1.
\end{equation}
\end{lemma}

\begin{proof}
Since $s \in \left] \tau_{m-1},\tau_m\right]$ and
$t_{\ell,\bi} \leq \tau_m$ implies $t_{\ell-1,\bi} \leq s$, we obtain
\begin{align*}
\int_0^1 E(\Yh_\bj^2(s)) \, ds
&\leq \scp{\xi}{h_\bj}^2 \cdot \int_0^1 \Gamma_\bj^2(s) \, ds
+
\sum_{\bi \in \I} \lambda_\bi/n_\bi \cdot 
\left( \sum_{\ell=0}^{n_\bi-1} a_{\bi,\bj}(t_{\ell,\bi})
\cdot \int_{t_{\ell,\bi}}^1 
\frac{\Gamma_\bj^2(s)}{\Gamma_\bj^2(t_{\ell,\bi})} \, ds
\right) \\
& \preceq
1/\mu_\bj \cdot
\left(\scp{\xi}{h_\bj}^2 + 
\sum_{\bi \in \I} \lambda_\bi/n_\bi \cdot 
\sum_{\ell=0}^{n_\bi-1} a_{\bi,\bj}(t_{\ell,\bi})
\right)
\end{align*}
from \eqref{g23} and Lemma \ref{l7}. It follows that
\[
\int_0^1 \ph(s) \, ds
\preceq
\|\xi\|^2 + 
\sum_{\bi \in \I} \lambda_\bi/n_\bi \cdot 
\sum_{\ell=0}^{n_\bi-1} (1 + E(\|\Xh(t_{\ell,\bi}) \|^2) ),
\]
see \eqref{g26}.
Use Lemma \ref{l6}  to complete the proof of \eqref{g22}.

Assume that $s < t$ with $s \in [\tau_{m-1},\tau_m]$ and
$t \in \left]\tau_{\kappa-1},\tau_\kappa\right]$ for $m \leq \kappa$.
Then
\begin{multline*}
E( \Yh_\bj (s) - \Yh_\bj(t))^2 \\
=
\left( 1 - \frac{\Gamma_\bj(t)}{\Gamma_\bj(s)} \right)^2 \cdot
E(\Yh_\bj^2(s)) 
+
\sum_{\bi \in \I} \lambda_\bi / n_\bi \cdot
\left(
\sum_{\ell \in \K_\bi(s,t)} 
\at_{\bi,\bj} (t_{\ell-1,\bi}) \cdot
\frac{\Gamma^2_\bj (t)}{\Gamma^2_\bj(t_{\ell-1,\bi})}
\right),
\end{multline*}
where
\[
\K_\bi (s,t) = 
\{ \ell \in \{1,\dots,n_\bi\} : t_{\ell,\bi} \in \left]\tau_m,\tau_\kappa
\right] \}
\]
if $s > \tau_{m-1}$ and
\[
\K_\bi (s,t) = 
\{ \ell \in \{1,\dots,n_\bi\} : t_{\ell,\bi} \in \left[\tau_m,\tau_\kappa
\right] \}
\]
if $s = \tau_{m-1}$. 
By Lemma \ref{l7}
\[
E( \Yh_\bj (s) - \Yh_\bj(t))^2
\preceq
\mu_\bj \cdot (t-s) \cdot E(\Yh_\bj^2(s))
+
\sum_{\bi \in \I} \lambda_\bi / n_\bi \cdot
\left(
\sum_{\ell \in \K_\bi(s,t)} 
\at_{\bi,\bj} (t_{\ell-1,\bi}) 
\right).
\]
Note that $\# \K_\bi(s,t) \leq 1 + n_\bi \cdot (t-s)$,
and apply \eqref{g26}
together with Lemma \ref{l6} to obtain
\begin{align*}
E\| \Xh(s) - \Xh(t)\|^2 
&\preceq
(t-s) \cdot \ph (s) +
\sum_{\bi \in \I} \lambda_\bi / n_\bi \cdot
\# \K_\bi(s,t)\\
&\preceq
(t-s) \cdot (1 + \ph (s)) + \sum_{\bi \in \I} \lambda_\bi / n_\bi,
\end{align*}
as claimed.
\end{proof}

In view of Lemma \ref{l6} and Lemma \ref{l8} we may
proceed as in the proof of Lemma \ref{l4}
to obtain the following error bound for piecewise constant 
interpolation of $\Xh$.

\begin{lemma}\label{l9}
For $\I, \J \subset \N^d$, $n \in \N^\I$, and $\bi \in \I$
\[
\sum_{\ell=1}^{n_\bi} \int_{t_{\ell-1,\bi}}^{t_{\ell,\bi}} 
E \|\Xh(t) - \Xh(t_{\ell-1,\bi})\|^2 \, dt 
\preceq 1/n_\bi + \sum_{\bi^\prime \in \I} \lambda_{\bi^\prime} /
n_{\bi^\prime}.
\]
\end{lemma}

\begin{proof}[Proof of Theorem \ref{t1}]
Recall the particular choice \eqref{g15} of the sets $\I$ and $\J$
and consider the corresponding It\^o-Galerkin approximation $\Xb$.
Because of Proposition \ref{l4} it suffices to show that
\begin{equation}\label{g21}
\int_0^1 E \|\Xb(t) - \Xh(t)\|^2 \, dt 
\preceq  \sum_{|\bi|_2 \leq I} \lambda_\bi / n_\bi. 
\end{equation}

For $\nu=1,2,3$ we define
\[
U^{(\nu)}_\bj(t) = \sum_{|\bi|_2 \leq I} \lambda_\bi^{1/2}
\cdot \int_0^t
 \sum_{\ell=0}^{n_\bi-1} V^{(\nu)}_{\bi,\bj,\ell}(s,t)
\cdot 1_{\left]t_{\ell,\bi},t_{\ell+1,\bi}\right]}(s) \, d\beta_\bi (s)
\]
with
\begin{align*}
V^{(1)}_{\bi,\bj,\ell} (s,t) 
&= \exp(-\mu_\bj \cdot(t-s)) \cdot
\scp{T_g(\Xb(s)) - T_g(\Xb(t_{\ell,\bi})) \cdot h_\bi}{h_\bj},\\
V^{(2)}_{\bi,\bj,\ell} (s,t) 
&= \exp(-\mu_\bj \cdot(t-s)) \cdot
\scp{T_g(\Xb(t_{\ell,\bi})) - T_g(\Xh(t_{\ell,\bi})) \cdot 
    h_\bi}{h_\bj},\\
V^{(3)}_{\bi,\bj,\ell} (s,t) 
&= \left( \exp(-\mu_\bj \cdot(t-s)) -
\frac{\Gamma_\bj(t)}{\Gamma_\bj(t_{\ell,\bi})} \right) \cdot
\scp{T_g(\Xh(t_{\ell,\bi})) \cdot h_\bi}{h_\bj}.
\end{align*}
Furthermore, we put
\[
U_\bj^{(4)}(t) =
\sum_{\bi \in \I \setminus \K_m}
\lambda_\bi^{1/2} \cdot
\frac{\Gamma_\bj(t)}{\Gamma_\bj(s_{m,\bi})} \cdot
\scp{T_g(\Xh(s_{m,\bi})) \cdot h_\bi}{h_\bj} 
\cdot (\beta_\bi(t) - \beta_\bi(s_{m,\bi}))
\]
and
\[
U_\bj^{(5)}(t) =
\sum_{\bi \in \K_m}
\lambda_\bi^{1/2} \cdot
\frac{\Gamma_\bj(t)}{\Gamma_\bj(s_{m,\bi})} \cdot
\scp{T_g(\Xh(s_{m,\bi})) \cdot h_\bi}{h_\bj} 
\cdot (\beta_\bi(\tau_m) - \beta_\bi(t))
\]
if $t \in \left]\tau_{m-1},\tau_m\right]$.
Then, by definition,
\begin{multline*}
\Yb_\bj (t) - \Yh_\bj (t)\\
=
\left( \exp (-\mu_\bj t) - \Gamma_\bj(t) \right) \cdot \scp{\xi}{h_\bj}
+ U_\bj^{(1)}(t) + U_\bj^{(2)}(t) + U_\bj^{(3)}(t) + 
U_\bj^{(4)}(t) - U_\bj^{(5)}(t).
\end{multline*}
We separately estimate the terms from the right-hand side of 
this equation.

Lemma \ref{l7} yields
\begin{equation}\label{g25}
\sum_{|\bj|_2 \leq J} \left( \scp{\xi}{h_\bj}^2 
\cdot \int_0^1 ( \exp(-\mu_\bj t) - \Gamma_\bj(t))^2 \, dt \right)
\preceq 1/n^* \preceq \sum_{|\bi|_2 \leq I} \lambda_\bi/n_\bi.
\end{equation}

By \eqref{g13} and Lemma \ref{l5}
\begin{equation}\label{g27}
\sum_{|\bj|_2 \leq J} E(U^{(1)}_\bj(t))^2 
\preceq
\sum_{|\bi|_2 \leq I}
\lambda_\bi \cdot \left(
\sum_{\ell=0}^{n_\bi-1}
\int_{t_{\ell,\bi}}^{t_{\ell+1,\bi}} 
E \| \Xb(s) -\Xb(t_{\ell,\bi})\|^2 \, ds \right)
\preceq
\sum_{|\bi|_2 \leq I}
\lambda_\bi / n_\bi.
\end{equation}

Put 
\[
f(s) = E\| \Xb(s) -\Xh(s)\|^2,
\]
which is finite because of \eqref{g14} and Lemma \ref{l6}.
By \eqref{g13}, Lemma \ref{l5}, and Lemma \ref{l9}
\begin{align}\label{g28}
&\sum_{|\bj|_2 \leq J} E(U^{(2)}_\bj(t))^2 \\
&\qquad \preceq
\sum_{|\bi|_2 \leq I}
\lambda_\bi \cdot \left( 
\sum_{\ell=0}^{n_\bi-1} 
\int_{t \wedge t_{\ell,\bi}}^{t \wedge t_{\ell+1,\bi}}
\left(
E\| \Xb(s) -\Xb(t_{\ell,\bi})\|^2 + E\| \Xh(s) -\Xh(t_{\ell,\bi})\|^2 + f(s)
\right)
\, ds \right)
\notag \\
&\qquad \preceq
\sum_{|\bi|_2 \leq I} \lambda_\bi / n_\bi + \int_0^t f(s)\, ds.
\notag
\end{align}

Suppose that $s \in \left]t_{\ell,\bi},t_{\ell+1,\bi}\right]$.
Then
\[
\left| \exp (-\mu_\bj(t-s)) - \exp(-\mu_\bj(t - t_{\ell,\bi})) \right|
\leq \exp (-\mu_\bj(t-s))\cdot \mu_\bj/n_\bi
\]
and therefore
\begin{align*}
&\int_s^1
\left( \exp(-\mu_\bj \cdot(t-s)) -
\frac{\Gamma_\bj(t)}{\Gamma_\bj(t_{\ell,\bi})} \right)^2 \, dt\\
&\qquad \preceq 
1/n_\bi +
\int_{t_{\ell,\bi}}^1
\left( \exp(-\mu_\bj \cdot(t-t_{\ell,\bi})) -
\frac{\Gamma_\bj(t)}{\Gamma_\bj(t_{\ell,\bi})} \right)^2 \, dt\\
&\qquad \preceq 1/n_\bi
\end{align*}
follows from Lemma \ref{l7}.
Hereby
\begin{align*}
& \int_0^1 E(U^{(3)}_\bj(t))^2 \, dt\\
&\preceq
\sum_{|\bi|_2 \leq I} \lambda_\bi \cdot
\left(
\int_0^1\!\int_0^t
\sum_{\ell=0}^{n_\bi-1}
\left( \exp(-\mu_\bj \cdot(t-s)) -
\frac{\Gamma_\bj(t)}{\Gamma_\bj(t_{\ell,\bi})} \right)^2
\cdot a_{\bi,\bj}(t_{\ell,\bi}) \cdot
1_{\left]t_{\ell,\bi},t_{\ell+1,\bi} \right]} (s)\, ds \, dt
\right)\\
&\preceq
\sum_{|\bi|_2 \leq I} \lambda_\bi / n_\bi  \cdot
\left(
\int_0^1\ \sum_{\ell=0}^{n_\bi-1}
a_{\bi,\bj}(t_{\ell,\bi}) \cdot
1_{\left]t_{\ell,\bi},t_{\ell+1,\bi} \right]} (s)\, ds
\right),
\end{align*}
which implies
\begin{equation}\label{g29}
\sum_{|\bj|_2 \leq J} \int_0^1 E(U^{(3)}_\bj(t))^2 \, dt
\preceq
\sum_{|\bi|_2 \leq I} \lambda_\bi / n_\bi, 
\end{equation}
see \eqref{g26} and Lemma \ref{l6}.

By the same facts,
\begin{equation}\label{g30}
\sum_{|\bj|_2 \leq J} 
E(U^{(4)}_\bj(t))^2 
\preceq
\sum_{|\bi|_2 \leq I} \lambda_\bi / n_\bi
\end{equation}
and
\begin{equation}\label{g31}
\sum_{|\bj|_2 \leq J} 
E(U^{(5)}_\bj(t))^2 
\preceq
\sum_{|\bi|_2 \leq I} \lambda_\bi / n_\bi, 
\end{equation}
if $t \in \left]\tau_{m-1},\tau_m\right]$.

Combining \eqref{g25}--\eqref{g31}
we obtain
\[
\int_0^r f(t) \, dt
\preceq
\sum_{|\bi|_2 \leq I} \lambda_\bi / n_\bi
+
\int_0^r \! \int_0^t f(s)\, ds \, dt.
\]
Finally, apply Gronwall's Lemma to derive
$\int_0^1 f(t)\, dt
\preceq \sum_{|\bi|_2 \leq I} \lambda_\bi / n_\bi$,
as claimed in \eqref{g21}.
\end{proof}

\section*{Acknowledgements}
\noindent
This work is supported by the Deutsche Forschungsgemeinschaft (DFG).

\section*{References}

{\small

\noindent
Appell, J., and Zabrejko, P.P. (1990),
Nonlinear Superposition Operators,
Cambridge Univ. Press, Cambridge.
\medskip

\noindent
Cambanis, S., and Hu, Y. (1996), 
Exact convergence rate of the Euler-Maruyama 
scheme, with application to sampling design, 
\SSR {\bf 59}, 211--240.
\medskip

\noindent
Da Prato, G., and Zabczyk, J. (1992),
Stochastic Equations in Infinite Dimensions,
Cambridge Univ. Press, Cambridge.
\medskip

\noindent
Davie, A. M., Gaines, J. (2001),
Convergence of numerical schemes for the solution of parabolic
partial differential equations, 
\MC {\bf 70}, 121--134.
\medskip

\noindent
Gy\"ongy, I. (1999),
Lattice approximations for stochastic quasi-linear parabolic partial
differential equations driven by space-time white noise II,
\PA {\bf 11}, 1--37.
\medskip

\noindent
Hausenblas, E. (2002),
Numerical analysis of semilinear stochastic evolution equations
in Banach spaces,
\JCAM {\bf 147}, 485--516.
\medskip

\noindent
Hausenblas, E. (2003),
Approximation for semilinear stochastic evolution equations,
\PA {\bf 18}, 141--186.
\medskip

\noindent
Hofmann, N., M\"uller-Gronbach, T., and Ritter, K. (2001),
The optimal discretization of stochastic differential equations.
\JC {\bf 17}, 117--153.
\medskip

\noindent
Kloeden, P., and Shott, S. (2001),
Linear-implicit strong schemes for It\^o-Galerkin approximations
of stochastic PDEs,
\JAMSA {\bf 14}, 47--53.
\medskip

\noindent
Manthey, R., and Zausinger, T. (1999),
Stochastic evolution equations in $L_\rho^{2\nu}$,
\SSR {\bf 66}, 37--85.
\medskip

\noindent
M\"uller-Gronbach, T. (2002),
Optimal uniform approximation of systems of stochastic differential 
equations, \AAP {\bf 12}, 664--690.
\medskip

\noindent
M\"uller-Gronbach, T. (2004),
Optimal pointwise approximation of SDEs
based on Brownian motion at discrete points,
\AAP {\bf 14}, 1605--1642.
\medskip

\noindent
M\"uller-Gronbach, T., and Ritter, K. (2006),
Lower bounds and nonuniform time discretization
for approximation of stochastic heat equations.
To appear in \FOCM
\medskip

\noindent
Millet, A., and Morien, P.-L. (2005),
On implicit and explicit discretization schemes for
parabolic SPDEs in any dimension,
\SPA {\bf 115}, 1073--1106.
\medskip

\noindent
Walsh, J. B. (2005),
Finite element methods for parabolic stochastic PDE's,
\PA {\bf 23}, 1--43.
\medskip

\noindent
Yan, Y. (2005),
Galerkin finite element methods for stochastic parabolic
partial differential equations,
\SJNA {\bf 43}, 1363--1384.
\medskip

}

\end{document}